\documentclass[12pt,a4paper]{article}
\usepackage{textcomp}
\usepackage{listings}
\usepackage[space=true]{accsupp}
\newcommand{\copyablespace}{\BeginAccSupp{method=hex,unicode,ActualText=00A0}\EndAccSupp{}}
\usepackage[T1]{fontenc}
\usepackage[english]{babel}
\usepackage[colorlinks = true,
            linkcolor = blue,
            urlcolor  = blue,
            citecolor = blue,
            anchorcolor = blue]{hyperref}
\usepackage{graphicx}
\usepackage{amsmath}
\usepackage{amssymb}
\usepackage{mathtools}
\usepackage[table,dvipsnames]{xcolor}
\usepackage{fancyhdr}
\usepackage{longtable}
\usepackage{multicol}
\usepackage{enumerate}
\usepackage{enumitem}
\setlist[itemize]{leftmargin=5.5mm}
\usepackage{multirow}
\usepackage{natbib}
\usepackage{caption}
\usepackage{subcaption}
\usepackage{xurl}
\usepackage[utf8]{inputenc}
\usepackage{tikz}
\usepackage{pgfplots}
\usepgfplotslibrary{colorbrewer}
\pgfplotsset{compat=newest, cycle list/Set1-8}
\usetikzlibrary {shapes.geometric}

\usetikzlibrary{pgfplots.statistics, pgfplots.colorbrewer}
\usepackage{pgfplotstable}

\pgfplotsset{
    box plot/.style={
        /pgfplots/.cd,
        black,
        only marks,
        mark=-,
        mark size=\pgfkeysvalueof{/pgfplots/box plot width},
        /pgfplots/error bars/y dir=plus,
        /pgfplots/error bars/y explicit,
        /pgfplots/table/x index=\pgfkeysvalueof{/pgfplots/box plot x index},
    },
    box plot box/.style={
        /pgfplots/error bars/draw error bar/.code 2 args={%
            \draw  ##1 -- ++(\pgfkeysvalueof{/pgfplots/box plot width},0pt) |- ##2 -- ++(-\pgfkeysvalueof{/pgfplots/box plot width},0pt) |- ##1 -- cycle;
        },
        /pgfplots/table/.cd,
        y index=\pgfkeysvalueof{/pgfplots/box plot box top index},
        y error expr={
            \thisrowno{\pgfkeysvalueof{/pgfplots/box plot box bottom index}}
            - \thisrowno{\pgfkeysvalueof{/pgfplots/box plot box top index}}
        },
        /pgfplots/box plot
    },
    box plot top whisker/.style={
        /pgfplots/error bars/draw error bar/.code 2 args={%
            \pgfkeysgetvalue{/pgfplots/error bars/error mark}%
            {\pgfplotserrorbarsmark}%
            \pgfkeysgetvalue{/pgfplots/error bars/error mark options}%
            {\pgfplotserrorbarsmarkopts}%
            \path ##1 -- ##2;
        },
        /pgfplots/table/.cd,
        y index=\pgfkeysvalueof{/pgfplots/box plot whisker top index},
        y error expr={
            \thisrowno{\pgfkeysvalueof{/pgfplots/box plot box top index}}
            - \thisrowno{\pgfkeysvalueof{/pgfplots/box plot whisker top index}}
        },
        /pgfplots/box plot
    },
    box plot bottom whisker/.style={
        /pgfplots/error bars/draw error bar/.code 2 args={%
            \pgfkeysgetvalue{/pgfplots/error bars/error mark}%
            {\pgfplotserrorbarsmark}%
            \pgfkeysgetvalue{/pgfplots/error bars/error mark options}%
            {\pgfplotserrorbarsmarkopts}%
            \path ##1 -- ##2;
        },
        /pgfplots/table/.cd,
        y index=\pgfkeysvalueof{/pgfplots/box plot whisker bottom index},
        y error expr={
            \thisrowno{\pgfkeysvalueof{/pgfplots/box plot box bottom index}}
            - \thisrowno{\pgfkeysvalueof{/pgfplots/box plot whisker bottom index}}
        },
        /pgfplots/box plot
    },
    box plot median/.style={
        /pgfplots/box plot,
        /pgfplots/table/y index=\pgfkeysvalueof{/pgfplots/box plot median index}
    },
    box plot width/.initial=1em,
    box plot x index/.initial=0,
    box plot median index/.initial=1,
    box plot box top index/.initial=2,
    box plot box bottom index/.initial=3,
    box plot whisker top index/.initial=4,
    box plot whisker bottom index/.initial=5,
}

\usepackage{framed}
\usepackage{booktabs}
\usepackage{caption}
\usepackage{float}
\usepackage{titlesec}
\usepackage{capt-of}
\definecolor{darkgreen}{rgb}{0.1, 0.5, 0.2}
\usepackage{array}
\usepackage{arydshln}
\newtheorem{definition}{Definition}

\definecolor{gray2}{rgb}{0.6,0.6,0.6}
\definecolor{lightgray2}{rgb}{0.8,0.8,0.8}
\title{An experimental approach: The graph of graphs} 
\author{Zsombor Sz\'adoczki$^{1,2,*}$, S\'andor Boz\'oki$^{1,2}$, L\'aszl\'o Sipos$^{3,4}$, Zs\'ofia Galambosi$^3$}
\date{}
\begin{document}
\pagenumbering{arabic}

\maketitle
\begin{center}
$*$ Corresponding author, 1111 Kende u. 13-17., Budapest, Hungary;\\ Email: szadoczki.zsombor@sztaki.hu\\
\bigskip
$^{1}$ Research Group of Operations Research and Decision Systems, \\
Research Laboratory on Engineering \& Management Intelligence \\
HUN-REN Institute for Computer Science and Control (SZTAKI), 1111 Kende u. 13-17., Budapest, Hungary;\\ Email: szadoczki.zsombor@sztaki.hu, bozoki.sandor@sztaki.hu\\
\bigskip
$^{2}$ Department of Operations Research and Actuarial Sciences \\
Corvinus University of Budapest, 1093 Fővám tér 8., Budapest, Hungary \\
\bigskip
$^3$ Hungarian University of Agriculture and Life Sciences,\\ Institute of Food Science and Technology, \\ Department of Postharvest, Commercial and Sensory Science, 1118, Villányi út 29-43., Budapest, Hungary;

\bigskip
$^4$ Institute of Economics,\\ HUN-REN Centre of Economic and Regional Studies, 1097 Tóth Kálmán utca 4., Budapest, Hungary;

\end{center}

\newpage
\renewcommand{\baselinestretch}{1.5} \normalsize

\begin{abstract}

\noindent
One of the essential issues in decision problems and preference modeling is the number of comparisons and their pattern to ask from the decision maker. We focus on the optimal patterns of pairwise comparisons and the sequence including the most (close to) optimal cases based on the results of a color selection experiment. In the test, six colors (red, green, blue, magenta, turquoise, yellow) were evaluated with pairwise comparisons as well as in a direct manner, on color-calibrated tablets in ISO standardized sensory test booths of a sensory laboratory. All the possible patterns of comparisons resulting in a connected representing graph were evaluated against the complete data based on 301 individual's pairwise comparison matrices (PCMs) using the logarithmic least squares weight calculation technique. It is shown that the empirical results, i.e., the empirical distributions of the elements of PCMs, are quite similar to the former simulated outcomes from the literature. The obtained empirically optimal patterns of comparisons were the best or the second best in the former simulations as well, while the sequence of comparisons that contains the most (close to) optimal patterns is exactly the same. In order to enhance the applicability of the results, besides the presentation of graph of graphs, and the representing graphs of the patterns that describe the proposed sequence of comparisons themselves, the recommendations are also detailed in a table format as well as in a Java application.

\end{abstract}


\renewcommand{\baselinestretch}{1.5} \normalsize

\section{Introduction}
\label{sec:1}

The use of pairwise (paired) comparisons is widespread in multicriteria decision making (MCDM) \citep{Greco2025}, preference measurement \citep{DavidsonFarquhar1976}, psychometry \citep{Thurstone1927}, sports \citep{Csato2021} as well as food science \citep{Sipos2025}.

In the practice of sensory testing, international standard methodologies are most commonly used to evaluate consumer preferences, which can be divided into three groups based on methodology (ranking, difference analysis, descriptive methods) \citep{ISO6658}. Among the standard methods, paired comparison tests (pairwise comparison, paired preference test) are one of the most common ranking techniques. The test is a forced choice test with two alternatives (two alternative forced choice, 2-AFC). In pairwise comparison, two products are evaluated by trained or expert assessors on the basis of the intensity of a sensory attribute, while in a paired preference test, consumers compare two products on the basis of their preference. By taking into account the dominant statistical characteristics - first-order error ($\alpha$), second-order error ($\beta$), proportion of discriminators (pd) - the difference or similarity of two products can be tested. Testing the difference between two products requires 30 assessors, while testing the similarity typically requires twice as many, 60 consumers, with equivalent sensitivity \citep{ISO5495}. The ranking of 3 to 6 products on the basis of one sensory attribute (intensity or liking) is aggregating the rankings and statistically evaluating them, where the comparison of the ranking sums of the products follows the logic of pairwise comparisons using non-parametric post hoc tests. In both test methods it should be noted that to determine the rankings of several properties, the sensory test should be carried out for each property separately \citep{ISO8587,Sipos2025}. In sensory practice, it is typically necessary to minimize the number of pairwise comparisons to prevent sensory fatigue and to maintain motivation, a principle that has been implemented in several standards \citep{ISO6658,ISO11136,ISO29842}.

The popular MCDM method, the Analytic Hierarchy Process (AHP) also applies pairwise comparisons (specifically pairwise comparison matrices (PCMs)) to evaluate the importance of different criteria, as well as to assess the performance of alternatives according to a given criterion \citep{Saaty1977,Saaty}. In the related methods, the lack of some comparisons, the incompleteness of data, is addressed in theory \citep{Ford1957,Harker1987}, and often occurs in practice as well \citep{BozokiCsatoTemesi,ISO29842}. The calculated outcome (the ranking of compared objects, items, alternatives, etc.) is highly dependent on the number of collected (known) comparisons, and on the arrangement of these comparisons. The latter refers to the pattern of the known comparisons, e.g., a pivotal alternative is compared to all the others, or each alternative is compared to the same number of other ones. This is often described by the (undirected) representing graph (graph of comparisons), where the vertices denote the alternatives, and there is an edge between two nodes if the comparison between the two appropriate elements is known \citep{Gass1998}.

Several special patterns of comparisons have been proposed in the literature, many of which include some sense of regularity of the representing graph \citep{WangTakahashi1998,Szadoczki2022}.

Taking a further step, the sequence of comparisons, connecting patterns with different numbers of known comparisons have also been examined \citep{FedrizziGiove2013}. This can be especially important when the decision maker is allowed to stop to answer the questions providing the appropriate pairwise comparisons. With an adequate sequence of comparisons, it still can be possible to estimate the decision maker's preferences sufficient enough. This problem is particularly interesting in the case of (large-scale) group decision making \citep{Tang2021}, when online questionnaires are used to gather the pairwise comparisons.

\cite{Szadoczki2025GraphofGraphs} analyzed all the possible patterns of comparisons for at most six alternatives, when the representing graph is connected. They determined the optimal patterns with the same number of known comparisons (and same number of alternatives), and proposed (partial) optimal sequences of comparisons based on the connection of the optimal patterns. The optimal patterns in their analysis were those that provided the closest weight vectors to the ones calculated from the complete data (complete PCMs) according to the Euclidean distance and the Kendall's $\tau$ measure. They applied graph of graphs to present their results, where every node of the analyzed graph is a graph itself, namely the representing graph of a given pattern of comparisons. There is an edge between two graphs in the graph of graphs if for each of the two represented patterns, the other one can be reached by the addition (deletion) of exactly one comparison (edge). Based on their results, and the findings of \cite{Gyarmati2023} and \cite{Szadoczki2023} the optimal filling in patterns seem to be robust for
\begin{itemize}
    \item the weight calculation technique,
    \item the level of inconsistency,
    \item the way of perturbation,
    \item the used distance metrics,
    \item the model using pairwise comparisons.
\end{itemize}

However, all of these studies applied simulations, although empirical pairwise comparison data can differ from simulated one significantly. We would like to fill in this research gap, and examine on empirical pairwise comparison matrices that how similar the optimal patterns of comparisons are to the simulated ones. We also use the concept of graph of graphs, and (to make it easier to follow) their components, from now on, are referred to using capital letters (e.g., GRAPH, EDGE, NODE, etc.) distinguishing them from the (representing) graphs (NODEs).

In this paper, the results of sensory-based experiments, namely a color-choice test is analyzed (for more details, see also \cite{Szadoczki2025}). The judgments of university students regarding six different colors were collected with pairwise comparisons using a four-item verbal category scale as well as direct scoring. It is examined that how close the results of the different patterns of comparisons are to the ones calculated from the complete PCMs.

The rest of the paper is organized as follows. Section~\ref{sec:2} contains the preliminaries related to pairwise comparisons. The methodology of the questionnaires and the evaluation of different patterns of comparisons is presented in Section~\ref{sec:3}. Section~\ref{sec:4} details the relation of optimal empirical and simulated patterns of comparisons. Section~\ref{sec:5} provides further discussion, and finally, Section~\ref{sec:6} concludes and raises research questions for the future.

\section{Pairwise comparisons and their analysis}
\label{sec:2}

Let us denote the number of alternatives (objects or items to be compared) by $n$. In this study, the focus is on the pairwise comparison matrix (PCM) application of paired comparisons. 

\begin{definition}[Pairwise comparison matrix (PCM)]
 The $n\times n$ matrix $\mathbf{A}=[a_{ij}]$ is a pairwise comparison matrix, if it is
 \begin{itemize}
     \item positive ($a_{ij}>0$ $\forall$  $i$,  $j$) and
     \item reciprocal ($1/a_{ij}  = a_{ji}$ $\forall$  $i$,  $j$).
 \end{itemize} 
\end{definition}

In the case of empirical problems, the decision makers usually provide inconsistent PCMs.

\begin{definition}[Consistent PCM]
A PCM is consistent if  $a_{ik}=a_{ij}a_{jk} \hspace{0.2cm} \forall i,j,k$. If a PCM is not consistent, then it is called inconsistent.
\end{definition}

There are a number of ways to measure the inconsistency of a PCM. Several inconsistency indices have been proposed in the literature \citep{Brunelli2018,Mazurek2023}, and many of their axiomatic properties have been investigated \citep{Brunelli2024} with special attention to their recommended thresholds \citep{Aguaron2003,Agoston2022}.

A variety of methods have been proposed to obtain a weight vector from a PCM that is the estimation of the decision maker's preference. In the case of a consistent PCM, the generally applied methods provide the same weight vector. However, for inconsistent empirical PCMs, the different methods can produce different weight vectors.

Some of the most popular weight calculation techniques are the right eigenvector method \citep{Saaty1977}, the componentwise reciprocal of the left eigenvector \citep{Johnson1976}, the least squares method \citep{Jensen1984,Bozoki2008b}, the technique based on the enumeration of spanning trees \citep{Tsyganok2010,Siraj2012,Lundy2017}, and the logarithmic least squares (geometric mean) method \citep{Crawford1985} (for further techniques and their comparison, see \cite{GolanyKress1993,ChooWedley2004,BajwaChooWedley2008,Dijkstra2013}). In this paper, the latter one, the logarithmic least squares (LLSM) method is applied.

\begin{definition}[Logarithmic Least Squares Method (LLSM)]
The weight vector $\mathbf{w}$ of an $n\times n$ PCM $\mathbf{A}$ determined by the LLSM is given by Equation~\ref{eq:1}.
\begin{equation}
\label{eq:1}
\min_{\mathbf{w}}    \sum_{i=1}^n\sum_{j=1}^n \left(\ln(a_{ij})-\ln\left(\frac{w_i}{w_j}\right)\right)^2 ,
\end{equation}
where $w_i$ is the $i$th coordinate of $\mathbf{w}$.
\end{definition}

In many applications the data is incomplete, i.e., some entries of the PCM are missing. That is called an incomplete pairwise comparison matrix (IPCM).

\begin{definition}[Incomplete pairwise comparison matrix (IPCM)]
An $n\times n$ matrix $\mathbf{A}=[a_{ij}]$ is an incomplete pairwise comparison matrix (IPCM) if:
\begin{itemize}
    \item $a_{ij}\in \mathbb{R}_+ \cup \{\ast\} \ \forall \ 1\leq i,j\leq n$,
    \item $a_{ji}=1/a_{ij}$ if $a_{ij}\in \mathbb{R}_+$,
    \item $a_{ji}=\ast$ if $a_{ij}=\ast$,
\end{itemize}
where $\ast$ denotes the missing elements, and $\mathbb{R}_+$ is the set of positive real numbers.
\end{definition}

The causes of the incompleteness can be various. One of them being the case that the decision maker has no time or willingness to provide all the comparisons, which can be especially interesting if the pairwise comparisons are provided via online questionnaires and the decision maker can simply stop answering the questions at any time \citep{Szadoczki2025GraphofGraphs}.

Many of the axiomatic properties of IPCMs can be studied suitably with the help of the representing graph (graph of comparisons) that shows the pattern of comparisons that is the main focus of the current paper.

\begin{definition}[Representing graph/Graph of comparisons]
An incomplete pairwise comparison matrix $\mathbf{A}$ can be represented by an undirected graph $G=(V,E)$, where:
\begin{itemize}
    \item the vertices $V=\{1,2,\ldots,n\}$ correspond to the alternatives,
    \item while the edge set $E$ represents the known elements of $\mathbf{A}$ outside the main diagonal:
    $$e_{ij} \in E \iff a_{ij}\neq\ast \ \text{and} \ i\neq j.$$
\end{itemize}
\end{definition}

The measurement of inconsistency \citep{Kulakowski2020}, and the most popular weight calculation techniques can be generalized to incomplete data as well \citep{BozokiTsyganok2019,Mazurek2022}. In the case of most techniques, the PCM is complemented minimizing an inconsistency index, and then the original method is applied \citep{Shiraishi1998,Shiraishi2002}. For the (generalized) LLSM, only the known elements are considered in Equation~\ref{eq:1}. However, it only provides a unique weight vector, if the related representing graph of the PCM is connected \citep{Bozoki2010}.

\section{Methodology}
\label{sec:3}

\subsection{Data: The conducted experiments}
\label{sec:3.1}

The empirical pairwise comparisons analyzed in this paper are obtained from a sensory-based color-choice test, also examined focusing on the applied optimal scale in \cite{Szadoczki2025}. The experiments were conducted in a laboratory that conforms to the International Standard for the design requirements of sensory tests  \citep{ISO8589}. They were carried out using color-calibrated tablets (Samsung Galaxy Tab A 2018, 10.5) with identical test geometry and maximum brightness settings in the sensory test booths.

In the tests, six colors ($n=6)$ were compared to each other in a pairwise manner using a four-item verbal category scale. Then the same six colors (red (R:189, G:62, B:57), green (R:90, G:151, B:90), blue (R:84, G:110, B: 183), magenta (R: 179, G: 55, B: 151), turquoise (R: 63, G: 185, B: 177), and yellow (R: 227, G: 203, B: 78) \citep{MylonasMacDonald2017}) were also directly evaluated by the participants on a 0-10 scale. The scale used to obtain numerical PCMs from the verbal results in this paper, is the one that is optimal according to \cite{Szadoczki2025}, i.e., the one that provides on average the smallest Euclidean distances between the weights calculated from the PCMs and the ones computed from direct scoring. See Table~\ref{tab:Scale} for details regarding the scale.

\begin{table}[ht!] \centering

\begin{tabular}{@{}lcc@{}}\toprule
	\multicolumn{1}{c}{Verbal expression} & Notation in optimization & Optimal numerical value \\ \hline
	Equally like them & $1$ & $1$ \\ 
	Like it a little bit more & $S$ & $1.5$  \\
	Like it moderately more & $M$ & $1.7$ \\
    Like it much more & $L$ & $2$ \\
	\bottomrule
	\end{tabular}
 \caption{The used scale for the conversion between verbal expressions and numerical values, see \cite{Szadoczki2025} for further details}
 \label{tab:Scale}
\end{table}

After certain filtering (that is related to the determination of the optimal scale), a total of 301 individuals preferences are contained in the examined data.

\subsection{Framework of the analysis of patterns of pairwise comparisons}
\label{sec:3.2}

To evaluate the different patterns of comparisons, the Euclidean distance ($d$) and the Kendall's $\tau$ measure are used:

\begin{equation}
\label{eq:2}
d(u,v)=\sqrt{\sum_{i=1}^{n}(u_i-v_i)^2},
\end{equation}

\begin{equation}
\label{eq:3}
\tau(u,v)=\frac{n_c(u,v)-n_d(u,v)}{n(n-1)/2},
\end{equation}

where $u$ and $v$ denote the weight vectors obtained from a certain pattern of comparisons and from the complete PCM, respectively. $u$ and $v$ are normalized by $  \sum_{i=1}^{n} u_i = 1$, and $  \sum_{i=1}^{n} v_i = 1$. $n_c(u,v)$ and $n_d(u,v)$ are the number of concordant and discordant pairs of the examined vectors, respectively. For the Euclidean distance, the smaller value is preferred, while for the Kendall's $\tau$, its higher value indicates a better performance of the given pattern.

The main assumptions applied in this paper are the same ones that are used in \cite{Szadoczki2025GraphofGraphs}:

\begin{enumerate}
    \item The comparisons can be chosen (they are not given a priori).
    \item An ‘optimal' pattern of comparisons (graph) is the one that provides the closest LLSM weight vector on average to the one calculated from the complete matrix according to the Euclidean distance and the Kendall's $\tau$ from the patterns with the same number of comparisons.
    \item \label{assumption:3} There is no prior information about the alternatives that should be compared, thus they can be handled in a symmetric way, only the non-isomorphic graphs are considered. This also means that it is assumed that the ‘reliability' and the weight of each comparison is the same.
\end{enumerate}

As a first step, because of Assumption~\ref{assumption:3}, a random permutation of the rows (and columns) of each empirical PCM was considered. This ensures that the alternatives cannot be distinguished, and the results can be compared to the former simulated outcomes. Otherwise it could be the case that a given color is much more popular than the others, and this would influence the results.

After the random permutation, the patterns of comparisons are considered by deleting the appropriate elements of the complete empirical PCMs. Then the LLSM weight vectors are calculated from all of those IPCMs (including the complete PCM as well), and the ones with the same number of comparisons (with the same number of edges in their representing graphs) are compared to each other based on the average Euclidean distance and the Kendall's $\tau$ measure calculated comparing them to the complete case. 

\section{Results}
\label{sec:4}

Based on the method described in Section~\ref{sec:3.2}, the results comparing the findings related to the empirical data, and the former simulated outcomes \citep{Szadoczki2025GraphofGraphs}, i.e., the two GRAPH of graphs for $n=6$, are presented in Figure~\ref{fig:GraphofGraphs}. Each NODE represent a given pattern of comparisons (a representing graph). The patterns with the same number of comparisons (edges, $e$) are presented in the same row. There is an EDGE between two NODEs if the two appropriate patterns of comparisons can be reached from each other by the addition/deletion of exactly one comparison (edge).

For a given number of comparisons ($e$), it often happens that the same pattern of comparisons turns out to be the best both according to the Euclidean distance, and the Kendall's $\tau$ measure. The NODEs representing these cases are colored green in Figure~\ref{fig:GraphofGraphs}. One can see that this is more common in the case of the simulated results. If there is an EDGE between two green NODEs (optimal patterns of comparisons), then that is also colored green (the two optimal patterns can be reached from each other by the addition/deletion of exactly one comparison). A detailed description of the empirically optimal graphs in comparison with the simulation results is presented in Table~\ref{tab:CompareOpt} that underlines the similarity of the outcomes for the two approaches.

\begin{table}[ht!] \centering
\footnotesize
\begin{tabular}{@{}cc@{}} \toprule
	 \multirow{ 2}{*}{\shortstack[c]{Number\\of comparisons}} &  \multirow{ 2}{*}{\shortstack[c]{Optimal empirical pattern compared to the simulated case}} \\
     &\\\hline
	\multirow{ 2}{*}{$e=5$} & \multirow{ 2}{*}{\shortstack[c]{The same optimal pattern (according to both measures):\\The star graph}} \\
    
    &\\
    
	\multirow{ 2}{*}{$e=6$} & \multirow{ 2}{*}{\shortstack[c]{The same optimal pattern (according to both measures):\\The $2$-regular $6$-cycle}}   \\

     &\\

	\multirow{ 2}{*}{$e=7$} & \multirow{ 2}{*}{\shortstack[c]{The same optimal pattern according to the Euclidean distance, a rank\\reversal between the best and the second best cases based on the Kendall's $\tau$}}\\

&\\
    
    \multirow{ 2}{*}{$e=8$} &  \multirow{ 2}{*}{\shortstack[c]{A rank reversal between the best and the second best\\patterns based on both measures}} \\

    &\\

     \multirow{ 2}{*}{$e=9$} & \multirow{ 2}{*}{\shortstack[c]{The same optimal pattern according to the Euclidean distance, a rank\\reversal between the best and the second best cases based on the Kendall's $\tau$}}\\

&\\

     \multirow{ 2}{*}{$e=10$} & \multirow{ 2}{*}{\shortstack[c]{The same optimal pattern according to the Euclidean distance, a rank\\reversal between the best and the second best cases based on the Kendall's $\tau$}}\\

&\\
         \multirow{ 2}{*}{$e=11$} & \multirow{ 2}{*}{\shortstack[c]{The optimal pattern (according to both measures)\\is the one that is optimal for the Euclidean distance in the simulations}}\\

&\\
        \multirow{ 2}{*}{$e=12$} & \multirow{ 2}{*}{\shortstack[c]{A rank reversal between the best and second best patterns based on both measures,\\but for each metric, a different pattern was the second best in the simulations}}\\

&\\
        
        $e=13$ & The same optimal pattern according to each of the two measures \\

	\bottomrule
	\end{tabular}
 \caption{The description of the optimal empirical patterns of comparisons according to the Euclidean distance and the Kendall's $\tau$ compared to the results of the simulations}
 \label{tab:CompareOpt}
\end{table}

A NODE is colored blue in Figure~\ref{fig:GraphofGraphs} when for the appropriate number of comparisons, either the Kendall’s $\tau$ and the Euclidean distance show different optimal cases, or it is the second best pattern according to both metrics, but the related optimal patterns cannot be reached from each other contrary to this one. In the case of the simulated GRAPH of graphs, for $e=12$, some of the patterns statistically provided the same results, thus three NODEs are highlighted by blue. For the empirical case, when $e=11$, there are different optimal patterns according to the two examined measures, but there is a third pattern that is the second best according to both measures, thus all three patterns are highlighted. EDGEs are colored blue in the case, when they are connecting blue NODES and they are part of the PATH that contains the most optimal (or second best) patterns from the $e=6$ case to the complete PCM ($e=n(n-1)/2=15$).

The spanning trees ($e=n-1=5$) are not considered in this PATH, as the optimal pattern in that case is the star graph, however, that is difficult to connect with any other optimal case with greater number of comparisons. On the other hand, the optimal pattern for $e=6$ can be reached with the addition of a single comparison from only one spanning tree, however, that is not among the best performers among the patterns with the same cardinality. At the same time, the preferences obtained from an IPCM with a spanning tree representing graph tend to be extremely unreliable. Thus, it is recommended to collect more comparisons whenever it is possible.

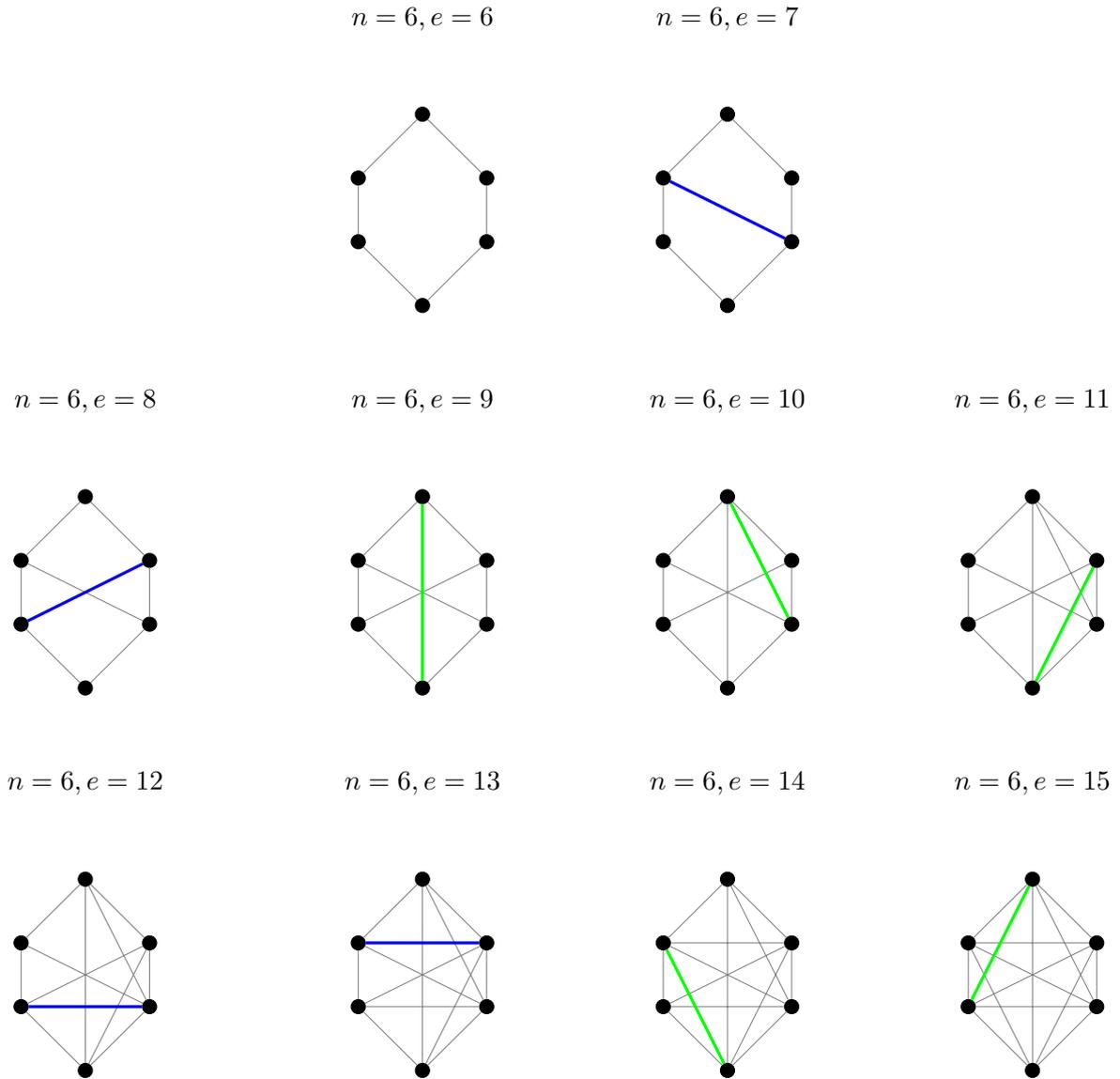
\begin{figure}[H]
\centering
\begin{tikzpicture}[every node/.style={circle,inner sep=2pt,draw=black,fill=black},scale=0.9]

\tikzstyle{block} = [circle,inner sep=2pt,draw=white,fill=white];


  \node (661) at (-4.75,-1.5) {};
  \node (662) at (-3.75,-2.5) {};
  \node (663) at (-5.75,-2.5) {};
  \node (664) at (-3.75,-3.5) {};
  \node (665) at (-5.75,-3.5) {};
  \node (666) at (-4.75,-4.5) {};
  
  \node [block] at (-4.75,0) (6) {\small $n=6, e=6$};
  
  \node (671) at (0,-1.5) {};
  \node (672) at (1,-2.5) {};
  \node (673) at (-1,-2.5) {};
  \node (674) at (1,-3.5) {};
  \node (675) at (-1,-3.5) {};
  \node (676) at (0,-4.5) {};
  
  \node [block] at (0,0) (7) {\small $n=6, e=7$};
  
  \node (681) at (-10,-7.5) {};
  \node (682) at (-9,-8.5) {};
  \node (683) at (-11,-8.5) {};
  \node (684) at (-9,-9.5) {};
  \node (685) at (-11,-9.5) {};
  \node (686) at (-10,-10.5) {};
  
  \node [block] at (-10,-6) (8) {\small $n=6, e=8$};
  
  \node (691) at (-4.75,-7.5) {};
  \node (692) at (-3.75,-8.5) {};
  \node (693) at (-5.75,-8.5) {};
  \node (694) at (-3.75,-9.5) {};
  \node (695) at (-5.75,-9.5) {};
  \node (696) at (-4.75,-10.5) {};
  
  \node [block] at (-4.75,-6) (9) {\small $n=6, e=9$};
  
  \node (6101) at (0,-7.5) {};
  \node (6102) at (-1,-8.5) {};
  \node (6103) at (1,-8.5) {};
  \node (6104) at (-1,-9.5) {};
  \node (6105) at (1,-9.5) {};
  \node (6106) at (0,-10.5) {};
  
  \node [block] at (0,-6) (10) {\small $n=6, e=10$};
  
  \node (6111) at (4.75,-7.5) {};
  \node (6112) at (3.75,-8.5) {};
  \node (6113) at (5.75,-8.5) {};
  \node (6114) at (3.75,-9.5) {};
  \node (6115) at (5.75,-9.5) {};
  \node (6116) at (4.75,-10.5) {};
  
  \node [block] at (4.75,-6) (11) {\small $n=6, e=11$};
  
  \node (6121) at (-10,-13.5) {};
  \node (6122) at (-9,-14.5) {};
  \node (6123) at (-11,-14.5) {};
  \node (6124) at (-9,-15.5) {};
  \node (6125) at (-11,-15.5) {};
  \node (6126) at (-10,-16.5) {};
  
  \node [block] at (-10,-12) (12) {\small $n=6, e=12$};
  
  \node (6131) at (-3.75,-14.5) {};
  \node (6132) at (-5.75,-14.5) {};
  \node (6133) at (-4.75,-13.5) {};
  \node (6134) at (-3.75,-15.5) {};
  \node (6135) at (-5.75,-15.5) {};
  \node (6136) at (-4.75,-16.5) {};
  
  \node [block] at (-4.75,-12) (13) {\small $n=6, e=13$};
  
  \node (6141) at (1,-14.5) {};
  \node (6142) at (-1,-14.5) {};
  \node (6143) at (0,-13.5) {};
  \node (6144) at (1,-15.5) {};
  \node (6145) at (-1,-15.5) {};
  \node (6146) at (0,-16.5) {};
  
  \node [block] at (0,-12) (9) {\small $n=6, e=14$};
  
   \node (6151) at (5.75,-14.5) {};
  \node (6152) at (3.75,-14.5) {};
  \node (6153) at (4.75,-13.5) {};
  \node (6154) at (5.75,-15.5) {};
  \node (6155) at (3.75,-15.5) {};
  \node (6156) at (4.75,-16.5) {};
  
  \node [block] at (4.75,-12) (9) {\small $n=6, e=15$};
  
  \draw [draw=black!50]
                        (661) -- (662)
                        (661) -- (663)
                        (662) -- (664)
                        (663) -- (665)
                        (664) -- (666)
                        (665) -- (666)
                        (671) -- (672)
                        (671) -- (673)
                        (672) -- (674)
                        (673) -- (675)
                        (675) -- (676)
                        (674) -- (676)
                        (673) -- (674)
                        (681) -- (682)
                        (681) -- (683)
                        (682) -- (684)
                        (683) -- (685)
                        (684) -- (686)
                        (685) -- (686)
                        (682) -- (685)
                        (683) -- (684)
                        (691) -- (692)
                        (691) -- (693)
                        (692) -- (694)
                        (693) -- (695)
                        (694) -- (696)
                        (695) -- (696)
                        (692) -- (695)
                        (693) -- (694)
                        (691) -- (696)
                        (6101) -- (6102)
                        (6101) -- (6103)
                        (6102) -- (6104)
                        (6103) -- (6105)
                        (6104) -- (6106)
                        (6105) -- (6106)
                        (6102) -- (6105)
                        (6103) -- (6104)
                        (6101) -- (6106)
                        (6101) -- (6105)
                        (6111) -- (6112)
                        (6111) -- (6113)
                        (6112) -- (6114)
                        (6113) -- (6115)
                        (6114) -- (6116)
                        (6115) -- (6116)
                        (6112) -- (6115)
                        (6113) -- (6114)
                        (6111) -- (6116)
                        (6111) -- (6115)
                        (6113) -- (6116)
                        (6121) -- (6122)
                        (6121) -- (6123)
                        (6122) -- (6124)
                        (6121) -- (6126)
                        (6123) -- (6125)
                        (6124) -- (6126)
                        (6125) -- (6126)
                        (6124) -- (6125)
                        (6121) -- (6124)
                        (6122) -- (6125)
                        (6123) -- (6124)
                        (6122) -- (6126)
                        (6131) -- (6132)
                        (6131) -- (6133)
                        (6132) -- (6134)
                        (6131) -- (6136)
                        (6133) -- (6134)
                        (6134) -- (6136)
                        (6135) -- (6136)
                        (6132) -- (6135)
                        (6131) -- (6135)
                        (6133) -- (6136)
                        (6131) -- (6134)
                        (6132) -- (6133)
                        (6134) -- (6135)
                        (6141) -- (6142)
                        (6141) -- (6143)
                        (6142) -- (6144)
                        (6143) -- (6144)
                        (6144) -- (6146)
                        (6145) -- (6146)
                        (6142) -- (6145)
                        (6141) -- (6145)
                        (6141) -- (6146)
                        (6141) -- (6144)
                        (6142) -- (6146)
                        (6142) -- (6143)
                        (6143) -- (6146)
                        (6144) -- (6145)
                        (6151) -- (6152)
                        (6151) -- (6153)
                        (6152) -- (6154)
                        (6153) -- (6155)
                        (6154) -- (6156)
                        (6155) -- (6156)
                        (6152) -- (6155)
                        (6153) -- (6154)
                        (6151) -- (6156)
                        (6151) -- (6154)
                        (6151) -- (6155)
                        (6152) -- (6156)
                        (6152) -- (6153)
                        (6153) -- (6156)
                        (6154) -- (6155);
\draw [very thick, draw=green]
                        (691) -- (696)
                        (6101) -- (6105)
                        (6113) -- (6116)
                        (6142) -- (6146)
                        (6153) -- (6155);
                        
\draw [very thick, draw=blue]
                        (673) -- (674)
                        (682) -- (685)
                        (6124) -- (6125)
                        (6131) -- (6132);

\end{tikzpicture}

\caption{The patterns of comparisons (graphs) related to the proposed sequence (the highlighted PATH in Figure~\ref{fig:GraphofGraphs}). The neighboring patterns can be reached from each other by the addition (deletion) of exactly one comparison (edge). The additional comparisons are highlighted in every step according to the color of the corresponding EDGE in the simulated part of Figure~\ref{fig:GraphofGraphs}.}
\label{fig:GraphAlongTheOptimalPATH}
\end{figure}

The representing graphs of the patterns of comparisons that are included in the optimal PATH of the GRAPH of graphs (the proposed filling in sequence of the PCM) are presented in Figure~\ref{fig:GraphAlongTheOptimalPATH}. This PATH is the same for the simulated and the empirical results as well. There are more possibilities in the case of the empirical results, i.e., one could choose another pattern for some of the $e$ values, but in that case this study proposes to use the one that is also included in the simulated results.

\newpage

\input{GraphofGraphs}

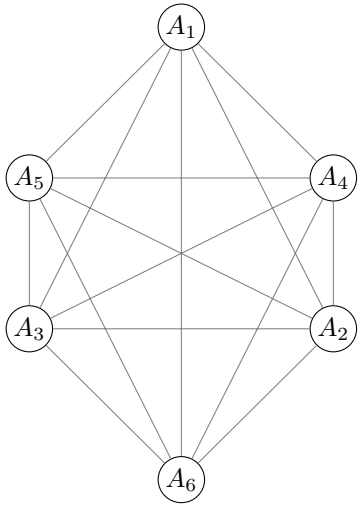
\begin{figure}[H]
\centering
\begin{tikzpicture}[every node/.style={circle,inner sep=1pt,draw=black,fill=white},scale=2]


   \node (6151) at (5.75,-14.5) {\footnotesize $A_4$};
  \node (6152) at (3.75,-14.5) {\footnotesize $A_5$};
  \node (6153) at (4.75,-13.5) {\footnotesize $A_1$};
  \node (6154) at (5.75,-15.5) {\footnotesize $A_2$};
  \node (6155) at (3.75,-15.5) {\footnotesize $A_3$};
  \node (6156) at (4.75,-16.5) {\footnotesize $A_6$};

  \draw [draw=black!50]
                        (6151) -- (6152)
                        (6151) -- (6153)
                        (6152) -- (6154)
                        (6153) -- (6155)
                        (6154) -- (6156)
                        (6155) -- (6156)
                        (6152) -- (6155)
                        (6153) -- (6154)
                        (6151) -- (6156)
                        (6151) -- (6154)
                        (6151) -- (6155)
                        (6152) -- (6156)
                        (6152) -- (6153)
                        (6153) -- (6156)
                        (6154) -- (6155);

\end{tikzpicture}

\caption{The complete graph with the labeling applied in Table~\ref{tab:FillingSequence}}
\label{fig:SampleLabel}
\end{figure}

The proposed filling in sequence is also described in Table~\ref{tab:FillingSequence} with the labeling presented in Figure~\ref{fig:SampleLabel} ($A_i$ stands for the $i$th alternative). Note that the appropriate filling sequence is equivalent to this one with the application of another labeling because of Assumption~\ref{assumption:3}. It is also worth mentioning that some of the comparisons are interchangeable, e.g., the first six comparisons can be made in any order (this is the starting point of the sequence), and the same is true for the last two comparisons (there is only one possible pattern for each case there, see Figure~\ref{fig:GraphofGraphs}). The ranks of the comparisons are colored according to the color of the corresponding additional edges in Figure~\ref{fig:GraphAlongTheOptimalPATH}. 

For instance, \#7 in Table~\ref{tab:FillingSequence} means that the comparison between $A_2$ and $A_5$ should be the seventh question in the questionnaire that collects the pairwise comparisons for a problem. The corresponding edge in Figure~\ref{fig:GraphAlongTheOptimalPATH} is the blue edge in the first row, and based on the labeling in Figure~\ref{fig:SampleLabel} that is between $A_2$ and $A_5$, indeed. The pattern obtained from the first seven comparisons this way is the one that is highlighted by a blue NODE in the $e=7$ row of the simulated results in Figure~\ref{fig:GraphofGraphs} from the 19 different possible patterns with seven pairwise comparisons. According to Table~\ref{tab:CompareOpt} this pattern provided the second best results in the simulations based on both measures, while it was the best according to the Kendall's $\tau$, and the second based regarding the Euclidean distance in the case of the empirical data.

\begin{table}[ht!]
\centering
\begin{tabular}{c||c|c|c|c|c|c|}

   &  $A_1$ \qquad & $A_2$   &  $A_3$   &  \ $A_4$ \  &    $A_5$    &    $A_6$    \\ \hline \hline
$A_1$  &     &  \color{green} \#10    &   \color{green} \#15   &  \#1  &  \#2   &  \color{green} \#9   \\ \hline
$A_2$  &     &      &  \color{blue}  \#12   &  \#3  &  \color{blue} \#7   &  \#4   \\ \hline
$A_3$  &     &      &      & \color{blue} \#8   &  \#5   &  \#6   \\ \hline
$A_4$  &     &      &      &        & \color{blue}   \#13     & \color{green}  \#11   \\ \hline
$A_5$  &     &      &      &        &         &  \color{green}  \#14      \\ \hline
$A_6$  &     &      &      &        &         &         \\     \hline
\end{tabular}
\caption{Proposed filling in sequence for $n = 6$, $e > 5$ based on the simulated and empirical results. The comparisons are colored according to the corresponding edges in Figure~\ref{fig:GraphAlongTheOptimalPATH}.} 
\label{tab:FillingSequence}
\end{table}

The results show that the assumed distributions in the simulations of \cite{Szadoczki2025GraphofGraphs} are not overwhelmingly different from the ones gained empirically in the color-choice experiment. The findings on the empirical data are quite similar to the simulated outcomes.

To further enhance the practical applicability of our proposed sequence, a Java application that provides assistance to use the results of the current paper as well as the results of \cite{Szadoczki2025GraphofGraphs} and \cite{Gyarmati2023} is available at \url{https://github.com/NNKDM8/Graph_of_graphs}. It displays the comparisons of the proposed patterns for given pairs of number of alternatives and number of comparisons ($n,e$), where the number of alternatives is between four and six.

\section{Discussion}
\label{sec:5}

When planning decision making tasks and experiments, it is necessary to consider the objective or subjective nature of the problem, the decision maker’s (assessor’s, judge’s) motivation, age, gender, cognitive abilities (concentration, memory, problem solving, understanding of test questions), level of proficiency in scale usage, the question sensitivity, the decision maker's involvement, etc. The main goal of our decision making experiment was to obtain the most information from the least test questions (to maximize the information retrieval), which is extremely important from both a decision theoretical and sensory perspective. Regarding the latter, determining the minimum number of test tasks is important in order to avoid sensory fatigue, reduce mental load, maintain attention and concentration levels. The modal and modal complexity of the test task largely determines the number of pairs and alternatives (products) to be tested. The degree and duration of sensory fatigue (decreased sensitivity) varies from sense to sense, starting with the least sensitive senses: hearing, vision, touch (somatosensory perception), taste, and smell.

The practical application of the results can be integrated into good sensory practices (GSP), and industrial product development becomes more efficient with the implementation of the fewest number of test tasks (less test tasks, less material, shorter development cycle, lower product development costs), thereby achieving more reliable results.

Although the examined empirical problem---the color-choice test---is highly subjective and specific regarding the number of alternatives ($n=6$), subjectivity is a common property of both decision making problems in general and sensory tests as well. The results support that the applied distribution of the elements of PCMs in the simulations of \cite{Szadoczki2025GraphofGraphs} is quite similar to those that occur in empirical problems. One could also argue that this strengthens the case for using the same simulation structure (that was first applied by \cite{Szadoczki2023}, and later adopted for instance by \cite{Csato2024}) related to pairwise comparisons independently from the examined question. The simulated results are examined for the cases when the number of alternatives is between four and six, and those types of problems are the most common in decision making. The empirical experiment was carried out for the largest (and regarding the possible patterns, the most interesting) case among these, while \cite{Szadoczki2025GraphofGraphs} argues that for larger problems the difference between certain patterns expected to be significantly smaller, which also confirms the relevance of the findings of this paper.

Similarly to any scientific work, the current study is not without limitations. It is difficult to quantify how similar or different results obtained from different type of problems can be, although the current results and the former simulations all point to the same direction. Most decision tasks do not apply Assumption~\ref{assumption:3}, and the labeled GRAPH of graphs could be examined as well. However, that means exponentially more patterns and, thus simulations can be computationally difficult, while for experiments this case seems to be even more problem-specific.

\section{Conclusion and further research}
\label{sec:6}

In this paper, an empirical experiment of a color selection task involving 301 individuals' pairwise and direct evaluation of six different colors was used to determine the optimal sequence of the comparisons, which was compared to previous runs on simulated data from the literature. The different patterns of comparisons with the same cardinality (same number of comparisons) were evaluated based on the Euclidean distance and the Kendall's $\tau$ measure between the logarithmic least squares weight vector obtained from them and from the complete data. The closer result to the complete data was considered to be better. It was assumed that there is no available prior information about the alternatives, thus, they were handled symmetrically (any case that only differs in the labeling of the alternatives considered to be the same pattern of comparisons).

It turned out that the empirical results---and so the empirical distributions of the elements of PCMs---are quite similar to the former simulated ones. The obtained empirically optimal patterns of comparisons were the best or the second best in the simulations as well, while the sequence of comparisons that contains the most (close to) optimal patterns is exactly the same.

The results were displayed using the concept of GRAPH of graphs (see Figure~\ref{fig:GraphofGraphs}), where the patterns of comparisons are represented by graphs (NODEs of the GRAPH), and each pair of graphs is connected if they can be reached from each other by the addition or deletion of exactly one comparison. In order to enhance the applicability of the results, besides each of the representing graphs of the patterns that describe the proposed sequence of comparisons (see Figure~\ref{fig:GraphAlongTheOptimalPATH}), the recommendations are also presented in a table format in Table~\ref{tab:FillingSequence} as well as in a Java application (\url{https://github.com/NNKDM8/Graph_of_graphs}).

Future research includes the investigation of the labeled GRAPH of graphs. Which patterns are the closest to the complete case, when there are some prior information about the possible performance of the alternatives? From a sensory perspective, it is worth investigating in what extent and how the proposed patterns of comparisons can be integrated into good sensory practices.

\section*{Acknowledgements}
The project identified by EKOP-CORVINUS-24-4-080 was realized with the support of the National Research, Development, and Innovation Fund provided by the Ministry of Culture and Innovation, as part of the University Research Scholarship Program announced for the 2024/2025 academic year. The research was supported by the National Research, Development and Innovation Office under Grants FK 145838 and TKP2021-NKTA-01 NRDIO.

\bibliographystyle{apalike} 
\bibliography{main}
\addcontentsline{toc}{section}{References}

\end{document}